\documentclass[reqno,centertags, 12pt]{amsart}

\usepackage{amssymb,amsmath,amsfonts,amssymb}%numbysec}
-\marginparwidth \oddsidemargin -9mm \evensidemargin \oddsidemargin
\usepackage{latexsym}
\advance\hoffset by 8mm

\def\@abssec#1{\vspace{.05in}\footnotesize \parindent .2in
{\bf #1. }\ignorespaces}
\newtheorem{theorem}{Theorem}[section]
\newtheorem{thm}[theorem]{Theorem}

\newcommand{\te}{\theta}
\newcommand{\J}{\mathcal{J}}
\newcommand{\Z}{\mathbb Z^2_+}

\def \Rm {\mathbb R}
\def \Tm {\mathbb T}

\renewcommand{\le}{\leqslant}
\renewcommand{\ge}{\geqslant}

\title{A simple energy pump for the surface quasi-geostrophic equation}
\author{Alexander Kiselev and Fedor Nazarov}
\thanks{Department of
Mathematics, University of Wisconsin, Madison, WI 53706, USA;
email: kiselev@math.wisc.edu, nazarov@math.wisc.edu}

\begin{document}
\begin{abstract}
We consider the question of growth of high order Sobolev norms of solutions of the conservative
surface quasi-geostrophic equation. We show that if $s>0$ is large then for every given $A$ there is exist small in $H^s$
initial data such that the corresponding solution's $H^s$ norm exceeds $A$ at some time. The idea of the
construction is quasilinear. We use a small perturbation of a stable shear flow. The shear flow can be shown to
create small scales in the perturbation part of the flow. The control is lost once the nonlinear effects become too large.
\end{abstract}

%{ \bf PRELIMINARY DRAFT}

\maketitle

\section{Introduction}\label{intro}

In this paper, we consider the surface quasi-geostrophic equation
\begin{equation}\label{sqg1}
\partial_t \theta = (u \cdot \nabla ) \theta, \,\,\,\theta(x,0)=\theta_0(x),
\end{equation}
$u = \nabla^\perp (-\Delta)^{-1/2}\theta,$ set on the torus $\Tm^2$ (which is equivalent to
working with periodic initial data in $\Rm^2$). Observe that the structure of the SQG equation is similar to the 2D Euler equation written for vorticity,
but the velocity is less regular in the SQG case ($u = \nabla^\perp (-\Delta)^{-1} \theta$ for the 2D Euler). The SQG equation comes from atmospheric science,
and can be derived via formal asymptotic expansion (assuming small Rossby and Ekman numbers) from a larger system of 3D Navier-Stokes
equations in a rotating frame coupled with temperature equation through gravity induced buoyancy force (see \cite{Pedl,Held}).
The equation \eqref{sqg1} describes evolution of the potential temperature on the surface, and its solution can be used to determine
the main order approximation for the solution of the full three dimensional problem.

%The setting that is natural from physical grounds is a sphere, in which case an extra term appears in \eqref{sqg1}.
%But this term is often not crucial for qualitative results, and for simplicity most papers consider \eqref{sqg1} on $\Rm^2$ with decaying or periodic
%boundary conditions.

In mathematical literature, the SQG equation was introduced for the first time by Constantin, Majda and Tabak in \cite{CMT}, where a parallel between the
structure of the conservative SQG equation and 3D Euler equation was drawn. Numerical experiments carried out in \cite{CMT} showed steep growth of the gradient
of solution in the saddle point scenario for the initial data, and suggested the possibility of singularity formation
in finite time. Subsequent numerical experiments \cite{OY} suggested the solutions stay regular. Later,
Cordoba \cite{Cord1} ruled out singularity formation in the scenario suggested by \cite{CMT}. Despite significant effort by many researchers,
whether blow up for the solutions of \eqref{sqg1} can happen in a finite time remains open. Moreover, there are no examples that exhibit just infinite growth in time
for some high order Sobolev norm. This paper is a step towards better understanding of of this phenomenon.

Before stating the main result, we would like to compare the situation with what is known for two-dimensional Euler equation, which in vorticity form
coincides with \eqref{sqg1} but the velocity is given by $u= \nabla^\perp (-\Delta)^{-1}\theta.$ The global existence of smooth solutions is known in this
case, and there is an upper bound on the gradient and higher order Sobolev norms of $\theta$ that is double exponential in time (see, e.g. \cite{MB}).
However the examples with actual growth are much weaker - the best current result is just superlinear in time (Denisov \cite{Den}, with earlier works
by Nadirashvili \cite{Nad} and Yudovich \cite{Yud2} giving linear or weaker rates of growth). It may appear that the SQG equation being more singular,
it should be easier to prove infinite growth in this case. However to prove infinite in time growth, one needs to produce an example of "stable instability", a controllable
mechanism of small scale production. This control is more difficult for the SQG than for two-dimensional Euler equation.

Let us denote $H^s$ the usual scale of Sobolev spaces on $\Tm^2.$
The main purpose of this short note is to show that the identically zero solution
is strongly unstable in $H^{s}$ for any $s$ sufficiently large. Namely, we will prove the following
\begin{thm}\label{th1}
Assume that $s$ is sufficiently large ($s \geq 11$ will do).
Given any $A>0,$ there exists $\theta_0$ such that $\|\theta_0\|_{H^s} \leq 1,$ but
the corresponding solution of \eqref{sqg1} satisfies
\begin{equation}\label{keyb} {\rm \limsup}_{t \rightarrow \infty} \|\theta(\cdot,t)\|_{H^s} \geq A. \end{equation}
\end{thm}
\it Remark. \rm 1. In our example, the initial data $\theta_0$ will be simply a trigonometric polynomial
with a few nonzero harmonics. Its size will be well controlled in any $H^s$ norm. \\
2. From the argument, it will be clear that it is not difficult to derive a lower bound on time
when the bound \eqref{keyb} is achieved. This time scales as a certain power of $A.$ \\
3. The arguments of Denisov \cite{Den} can also be used to produce similar result, with better control of constants and time - but in a different scenario.
Denisov considers perturbation of an explicitly given saddle point flow.
In this note, we will consider a technically simpler but less singular case of a shear flow.

\section{The Proof}\label{proof}

We shall view the solutions $\te(x,t)$ of \eqref{sqg1} as sequences of Fourier coefficients $\hat{\te}_k$, $k=(k_1,k_2)\in\mathbb Z^2$.
On the Fourier side, after symmetrization, our solution satisfies the following equation.
\begin{equation}\label{foursqg}
\frac{d}{dt}\hat{\te}_k = \frac12 \sum_{l+m=k}(l\wedge m)\left(\frac 1{|l|}-\frac 1{|m|}\right)\hat{\te}_l\hat{\te}_m, \,\,\,\hat{\te}_k(0) = (\hat{\te_0})_k.
\end{equation}
Here $l \wedge m = l_1 m_2 - l_2 m_1.$

Our initial data $\theta_0$ will be just a simple trigonometric polynomial $p,$ given by $\hat{p}_e=\hat{p}_{-e}=1$,
$\hat{p}_g=\hat{p}_{g+e}=\hat{p}_{-g}=\hat{p}_{-g-e}=\tau$ where $e=(1,0)$, $g=(0,2)$, and $\tau=\tau(A)>0$
is a small parameter to be chosen later. Then it follows from \eqref{foursqg} that the solution is an even real-valued function
with $\hat{\te}_0=0$ for all times. Moreover, $\hat{\te}_k(t) \equiv 0$ whenever $k_2$ is odd.

We have two easy to check conservation laws: $\sum_k \hat{\te}_k(t)^2=2+4\tau^2$ and
$\sum_k \frac{\hat{\te}_k^2(t)}{|k|}=2+2\tau^2\left(\frac 12+\frac 1{\sqrt 5}\right)$. 
After subtraction, we obtain that 
\[ \sum_{k} \hat{\te}_k(t)^2 \left(1- \frac{1}{|k|}\right) = \left(3 -\frac{2}{\sqrt{5}}\right)\tau^2, \]
for all $t \geq 0.$ Since $\hat{\te}_{\pm g/2}(t) =0,$ this implies 
\[ \sum_{k \ne \pm e} \hat{\te}_k(t)^2 \leq \left(\frac{\sqrt{2}}{\sqrt{2}-1}\right)\left(3 -\frac{2}{\sqrt{5}}\right)\tau^2 
\leq 10 \tau^2 \]
for all times. Then the first conservation law also implies $\hat{\te}_e(t) \in (1-8\tau^2,1+2\tau^2) \subset (1/2,2)$ for all times, 
provided that $\tau$ is sufficiently small. 

%From these laws, it follows
%that $\hat{\te}_e=\hat{\te}_{-e}\in \left(\frac 12,2\right)$ and $\sum_{k\neq\pm e}\hat{\te}_k^2\le 10\tau^2$ for all times.

Consider the quadratic form
$$
\J(\hat{\te})=\sum_{k\in \Z}\Phi(k)\hat{\te}_k\hat{\te}_{k+e}\,.
$$
We have
\begin{multline*}
\frac{d}{dt}\J(\hat{\te})=
\frac12 \sum_{k\in\Z}\Phi(k)\left[
\hat{\te}_k\sum_{l+m=x+e\,,\,l,m\ne\pm e}(l\wedge m)\left(\frac 1{|l|}-\frac 1{|m|}\right)\hat{\te}_l\hat{\te}_m
\right.
\\
+\left.
\hat{\te}_{k+e}\sum_{l+m=k\,,\,l,m\ne\pm e}(l\wedge m)\left(\frac 1{|l|}-\frac 1{|m|}\right)\hat{\te}_l\hat{\te}_m
\right]
+
\\
\frac12 \hat{\te}_e\sum_{k\in\Z}(e\wedge k)\Phi(k)\left[
\left(1-\frac1{|k|}\right)\hat{\te}_k^2-
\left(1-\frac1{|k+2e|}\right)\hat{\te}_k\hat{\te}_{k+2e}\right.
\\
-\left.
\left(1-\frac1{|k+e|}\right)\hat{\te}_{k+e}^2+
\left(1-\frac1{|k-e|}\right)\hat{\te}_{k+e}\hat{\te}_{k-e}
\right]
\equiv \sigma+\Sigma\,.
\end{multline*}
where $\sigma$ denotes the first sum and $\Sigma$ the second.
Since for $l+m=k$, we have $\left|(l\wedge m)\left(\frac 1{|l|}-\frac 1{|m|}\right)\right|\le 2|k|$ and $|k+e|\asymp |k|$ for $k\in \Z$, we conclude that
$$
|\sigma|\le \left(\sum_{k\in\Z}|k||\Phi(k)||\hat{\te}_k|\right)\left(\sum_{l\ne\pm e}\hat{\te}_l^2\right)
\le
C\tau^2\sum_{k\in\Z}|k||\Phi(k)||\hat{\te}_k|\,.
$$
%%
%Set $\Phi(k)=k_1+\frac12$.
%As we will see below, the sum in the expression for $\Sigma$ is positive.
%Therefore, since $\hat{\te}_e(t) \geq 1/2$ for all $t \geq 0,$ 
On the other hand, $\Sigma$ can be rewritten as  
\begin{multline*}
\hat{\te}_e\sum_{k_2>0}k_2\sum_{k_1\in\mathbb Z}\frac 14\times
\\
\left[
(\Phi(k-e)-\Phi(k-2e))\left(1-\frac1{|k-e|}\right)\hat{\te}_{k-e}^2+
(\Phi(k+e)-\Phi(k))\left(1-\frac1{|k+e|}\right)\hat{\te}_{k+e}^2
\right.
+
\\ \left.
2\left\{\Phi(k))\left(1-\frac1{|k-e|}\right)-
\Phi(k-e)\left(1-\frac1{|k+e|}\right)\right\}\hat{\te}_{k-e}\hat{\te}_{k+e}
\right]
\end{multline*}
Now let $\Phi(k) = k_1+\frac12.$ 
We get the sum of quadratic forms with the coefficients
$$
1-\frac{1}{\sqrt{(k_1-1)^2+k_2^2}}\ ,\ 1-\frac{1}{\sqrt{(k_1+1)^2+k_2^2}}
$$
at the squares
and
$$
\left(k_1+\frac12\right)\left(1-\frac{1}{\sqrt{(k_1-1)^2+k_2^2}}\right)-
\left(k_1-\frac12\right)\left(1-\frac{1}{\sqrt{(k_1+1)^2+k_2^2}}\right)
$$
at the double product.

A straightforward computation shows that when $k_1=0$, this form is degenerate and when $k_1\ne 0$, it is strictly positive definite
and dominates $\frac c{|k|^3}(\hat{\te}_{k-e}^2+\hat{\te}_{k+e}^2)$.

Using that fact that $\hat{\te}_e(t) \geq 1/2$ for all times, we obtain
$$
\Sigma\ge c\sum_{k\in\Z} \frac{\hat{\te}_k^2}{|k|^3}\,.
$$
Now there are several possibilities.

A) At some time $t,$ we will have $\sum_{k\in\Z}|k|^2|\hat{\te}_k|\ge \sum_{k\in\Z}|k||\Phi(k)||\hat{\te}_k|\ge \tau^{1/2}\,.$

Observe that
$$
\sum_{k\in\Z}|k|^2|\hat{\te}_k|\le\left(\sum_{k\in\Z}\hat{\te}_k^2\right)^{1/3}
\left(\sum_{k\in\Z}|k|^{21}\hat{\te}_k^2\right)^{1/6}\left(\sum_{k\in\Z}|k|^{-3}\right)^{1/2}
$$
Then, since $\sum_{k\in\Z}\hat{\te}_k^2\le 10\tau^2$, we get
that the $H^{11}$ norm of the solution gets large: $\|\theta\|_{H^{11}} \geq C \tau^{-1/6}.$

B) The case (A) never occurs but $\sum_{k\in\Z}|k|^{-3}\hat{\te}_k^2$ becomes comparable with
$\tau^{5/2}$. Note that until this moment $\J(\hat{\te})$ increases from its initial value about $\tau^2$. Also, $\J(\hat{\te})\le \sum_{k\in\Z}|k|\hat{\te}_k^2$.

Thus, in this case, we use
$$
\sum_{k\in\Z}|k|\hat{\te}_k^2\le\left(\sum_{k\in\Z}|k|^{-3}\hat{\te}_k^2\right)^{5/6}
\left(\sum_{k\in\Z}|k|^{21}\hat{\te}_k^2\right)^{1/6}\,.
$$
and, again, it follows that that the $H^{11}$ norm becomes large: $\|\theta\|_{H^{11}} \geq C \tau^{-\frac{1}{12}}$.

At last, if neither (A), nor (B) occur, then $\J(\hat{\te})$ grows without
bound and the $H^{1/2}$-norm gets large eventually. Now given $A$ just choose $\tau$ sufficiently small and Theorem~\ref{th1} follows. \\

\noindent {\bf Acknowledgement.} \rm Research of AK has been
supported in part by the NSF-DMS grant 0653813. Research of FN has
been partially supported by the NSF-DMS grant 0800243.

\end{document}